\documentclass{proc-l}

\usepackage{amssymb}

\theoremstyle{plain}
\newtheorem{theorem}{Theorem}[section]
\newtheorem{lemma}[theorem]{Lemma}

\newtheorem{milnorthmone}{Theorem 1 (Milnor)}

\newtheorem{milnorlemone}{Lemma 1 (Milnor)}

\newtheorem{milnorlemtwo}{Lemma 2 (Milnor)}

\newtheorem{milnorthmtwo}{Theorem 2 (Milnor)}

\theoremstyle{remark}
\newtheorem{remark}[theorem]{Remark}
\newtheorem*{acknowledgment}{Acknowledgment}

\newcommand{\ab}{\ensuremath{A \# B}}
\newcommand{\ds}{\ensuremath{D_1 \# ... \# D_s}}
\newcommand{\mspo}{\ensuremath{M_1 \# ... \# M_r}}
\newcommand{\snms}{\ensuremath{S^3 \# M_1 \# ... \# M_{r-1}}}

\newcommand{\g}{\ensuremath{\gamma}}
\newcommand{\w}{\ensuremath{\omega}}
\newcommand{\tw}{\ensuremath{\tilde \omega}}
\newcommand{\s}{\ensuremath{\sigma}}
\newcommand{\F}{\ensuremath{\Sigma}}

\newcommand{\ta}{\ensuremath{\tilde A}}
\newcommand{\tb}{\ensuremath{\tilde B}}
\newcommand{\tf}{\ensuremath{\tilde F}}
\newcommand{\tm}{\ensuremath{\tilde M}}
\newcommand{\tp}{\ensuremath{\tilde P}}
\newcommand{\ts}{\ensuremath{\tilde S}}
\newcommand{\tT}{\ensuremath{\tilde T}}
\newcommand{\tx}{\ensuremath{\tilde X}}

\newcommand{\brts}{\ensuremath{[\tilde S]}}

\newcommand{\va}{\ensuremath{V_A}}
\newcommand{\bva}{\ensuremath{\partial V_{A}}}
\newcommand{\vaj}{\ensuremath{V_{A,j}}}
\newcommand{\bvaj}{\ensuremath{\partial V_{A,j}}}
\newcommand{\tva}{\ensuremath{\tilde V_A}}
\newcommand{\vb}{\ensuremath{V_B}}
\newcommand{\bvb}{\ensuremath{\partial V_{B}}}
\newcommand{\vbj}{\ensuremath{V_{B,j}}}
\newcommand{\bvbj}{\ensuremath{\partial V_{B,j}}}
\newcommand{\tvb}{\ensuremath{\tilde V_B}}

\newcommand{\mxs}{\ensuremath{M \times S^1}}
\newcommand{\tmxs}{\ensuremath{\tilde M \times S^1}}
\newcommand{\sxs}{\ensuremath{S^2 \times S^1}}
\newcommand{\sxt}{\ensuremath{S^2 \times S^1 \times S^1}}

\newcommand{\boa}{\ensuremath{b_1(A)}}
\newcommand{\bota}{\ensuremath{b_1(\tilde A)}}
\newcommand{\botb}{\ensuremath{b_1(\tilde B)}}
\newcommand{\bom}{\ensuremath{b_1(M)}}
\newcommand{\botm}{\ensuremath{b_1(\tilde M)}}
\newcommand{\bop}{\ensuremath{b_1(P)}}

\newcommand{\btm}{\ensuremath{b_2(M)}}
\newcommand{\btx}{\ensuremath{b_2(X)}}

\newcommand{\bpx}{\ensuremath{b_2^+(X)}}
\newcommand{\bnx}{\ensuremath{b_2^-(X)}}

\newcommand{\hoa}{\ensuremath{H_1(A)}}
\newcommand{\hota}{\ensuremath{H_1(\tilde A)}}
\newcommand{\hotb}{\ensuremath{H_1(\tilde B)}}

\newcommand{\htx}{\ensuremath{H_2(X)}}

\newcommand{\homthos}{\ensuremath{H_1(M) \otimes H_1(S^1)}}
\newcommand{\htmthzs}{\ensuremath{H_2(M) \otimes H_0(S^1)}}

\newcommand{\tmsS}{\ensuremath{\tilde M \setminus \Sigma}}

\newcommand{\R}{\ensuremath{\mathbb R}}
\newcommand{\Z}{\ensuremath{\mathbb Z}}

\newcommand{\fund}{\ensuremath{\pi_1}}
\newcommand{\pii}{\ensuremath{\pi_i}}

\begin{document}

\title{On the Asphericity of a Symplectic $M^3 \times S^1$}

\author{John D. McCarthy}
\address{Department of Mathematics \\ Michigan State University \\ East Lansing, MI 48824}
\email{mccarthy@math.msu.edu}
\keywords{Symplectic, manifold, aspherical}
\subjclass{Primary 53C15; Secondary 57M50, 57N10, 57N13}
\date{April 6, 1999}
\commby{Ronald A. Fintushel}

\begin{abstract} C. H. Taubes asked whether a closed (i.e. compact and without boundary)
connected oriented three-dimensional manifold whose product with a circle admits a symplectic
structure must fiber over a circle. An affirmative answer to Taubes' question would imply that
any such manifold either is diffeomorphic to the product of a two-sphere with a circle or is
irreducible and aspherical. In this paper, we prove that this implication holds up to connect
sum with a manifold which admits no proper covering spaces with finite index. It is pointed out
that Thurston's geometrization conjecture and known results in the theory of three-dimensional
manifolds imply that such a manifold is a three-dimensional sphere. Hence, modulo the present
conjectural picture of three-dimensional manifolds, we have shown that the stated consequence
of an affirmative answer to Taubes' question holds.
\end{abstract}

\maketitle

\section{Introduction}
\label{sec:intro}

An interesting question in symplectic topology, which was posed by C. H. Taubes, concerns the
topology of closed (i.e. compact and without boundary) connected oriented three-dimensional
manifolds whose product with a circle admits a symplectic structure. The only known examples of
such manifolds are those which fiber over a circle. Taubes asked whether these examples are the
only examples of such manifolds.

Let $M$ be a closed oriented $3$-manifold which fibers over the circle
$S^1$. By definition, $M$ is the total space of a locally trivial fiber bundle whose base
space is the circle $S^1$ and whose fiber
$F$ is a closed oriented $2$-manifold. Let $\pi : M \rightarrow S^1$ be the projection of this
bundle. By replacing $\pi$ by a lift $\tilde{\pi} : M \rightarrow S^1$ of
$\pi$ to an appropriate covering space $S^1 \rightarrow S^1$ of $S^1$ with finite index, we may
assume that
$F$ is connected.

We may pull back the universal covering space $\R \rightarrow S^1$ of
$S^1$ through $\pi$ to obtain an infinite  cyclic covering space $Z$ of $M$ which fibers over
the line
$\R$ with fiber $F$. It follows that $Z$ is diffeomorphic to $F \times \R$. Thus, the
universal covering space $\tm$ of $M$ is diffeomorphic to $\tf \times \R$, where
$\tf$ is the universal covering space of $F$. As is well known, $\tf$ is diffeomorphic to
$S^2$, if $g = 0$, and to $\R^2$, if $g \geq 1$.

Thus, if $g = 0$, $M$ is diffeomorphic to the product $\sxs$ of $S^2$ with
$S^1$. (Note that every orientation preserving diffeomorphism of $S^2$ is isotopic to the
identity. Hence, there are exactly two $S^2$-bundles over $S^1$, the product $\sxs$ of $S^2$
with $S^1$ (i.e the mapping torus of the identity map of $S^2$) and the twisted product of
$S^2$ with $S^1$ (i.e the mapping torus of the antipodal map of $S^2$). Since $M$ is
orientable and the twisted product of
$S^2$ with $S^1$ is nonorientable, $M$ can not be diffeomorphic to a twisted product of $S^2$
with $S^1$.) On the other hand, if $g \geq 1$, then $\tm$ is diffeomorphic to $\R^3$ and,
hence, $M$ is irreducible (i.e. every (tame) $2$-sphere in $M$ bounds a $3$-ball) and
aspherical (i.e. the higher homotopy groups
$\pii(M), i > 1,$ of $M$ are zero).

Suppose now that $M$ is a closed connected oriented $3$-manifold whose product $\mxs$ with the
circle $S^1$ admits a symplectic structure. From the previous observations, it follows that an
affirmative answer to Taubes' question would imply that either $M$ is diffeomorphic to $\sxs$
or $M$ is irreducible and aspherical. Note that the universal cover of $\mxs$ is diffeomorphic
to
$\tm \times \R$. It follows that $\mxs$ is aspherical if and only if $M$ is aspherical.

In this paper, we shall prove the following result.

\begin{theorem} Let $M$ be a closed connected oriented $3$-manifold whose product $\mxs$ with
the circle $S^1$ admits a symplectic structure. Then $M$ has a unique connect sum
decomposition $\ab$ satisfying the following conditions: (i) the first Betti number $\boa$ of
$A$ is at least $1$, (ii) either $A$ is diffeomorphic to $\sxs$ or $A$ is irreducible and
aspherical, and (iii) every connected covering space of $B$ with finite index is trivial.
\label{theorem:mainthm} \end{theorem}

\begin{remark} Given that $\boa \geq 1$ and $A$ is irreducible, the statement that $A$ is
aspherical is redundant. Indeed, any compact oriented irreducible
$3$-manifold with an infinite fundamental group is aspherical, (see Chapter IV, Section 2 of
\cite{mc}). Since our main interest in this paper concerns asphericity, we include the
adjective ``aspherical'' in Theorem
\ref{theorem:mainthm} for emphasis. On the other hand, we include the adjective
``irreducible'' in Theorem
\ref{theorem:mainthm} in order to achieve the uniqueness clause of this theorem.
\label{remark:redundancy} \end{remark}

\begin{remark} Note that Theorem \ref{theorem:mainthm} addresses the question of the extent to
which a symplectic $\mxs$ must be aspherical. Furthermore, the previous observations imply
that an affirmative answer to Taubes' question would imply that $B$ is a
$3$-sphere.
\label{remark:sphere} \end{remark}

\begin{remark} Note that the factor $B$ in Theorem \ref{theorem:mainthm} is a homology sphere.
Indeed, a manifold is an homology sphere if and only if every connected cyclic covering space
of the manifold with finite index is trivial. In addition, if the fundamental group
$\fund(B)$ of this factor $B$ is finite, then $B$ must be a homotopy sphere.
\label{remark:homologysphere} \end{remark}

\begin{remark} We recall that a group $G$ is residually finite if every nontrivial element of
$G$ is mapped nontrivially to some finite quotient group of $G$. In particular, if
$N$ is a connected
$3$-manifold with a nontrivial residually finite fundamental group, then
$N$ has a nontrivial connected covering space with finite index. Hence, if the factor $B$ in
Theorem
\ref{theorem:mainthm} has a residually finite fundamental group, it must be a homotopy sphere.
\label{remark:resfinite} \end{remark}

\begin{remark} A theorem of Thurston's \cite{th} states that the fundamental group of a Haken
manifold is residually finite. As pointed out by Hempel \cite{h}, it follows that this theorem
extends to the class of all $3$-manifolds whose prime factors either are virtually Haken or
have finite or cyclic fundamental groups. As also pointed out by Hempel
\cite{h}, it is unsolved whether this class includes all closed 3-manifolds (cf. \cite{th},
section 6). On the other hand, Thurston's ``geometrization conjecture'' implies that every
closed 3-manifold lies in this class. Hence, Theorem
\ref{theorem:mainthm} and Thurston's ``geometrization conjecture'' imply that $B$ must be a
homotopy sphere. On the other hand, as is well-known, Thurston's ``geometrization conjecture''
also implies the Poincare conjecture. Hence, Theorem \ref{theorem:mainthm} and Thurston's
``geometrization conjecture'' imply that $B$ must be diffeomorphic to a sphere.  It follows
that Theorem \ref{theorem:mainthm} and Thurston's ``geometrization conjecture'' imply that
either
$M$ is diffeomorphic to $\sxs$ or
$M$ is aspherical and, hence, that either $\mxs$ is diffeomorphic to $\sxt$ or $\mxs$ is
aspherical.
\label{remark:asphericity} \end{remark}

\begin{remark} In \cite{k2}, Dieter Kotschick conjectured that a general
$4$-manifold which is symplectic for both choices of orientation must be either ruled or
aspherical. Note that if $\mxs$ is symplectic, then it is symplectic for both choices of
orientation. Furthermore, note that $\sxt$ is ruled (being an $S^2$ bundle over $S^1 \times
S^1$). We conclude, from the previous remark, that Theorem \ref{theorem:mainthm} and
Thurston's ``geometrization conjecture'' imply that Kotschick's conjecture holds for a general
$4$-manifold of the form $\mxs$.
\label{remark:kotschickconjecture} \end{remark}

Here is an outline of the paper. In Section 1, we shall prove some technical lemmas about
connect sum decompositions of closed connected oriented $3$-manifolds whose product with a
circle admits a symplectic structure. In section 2, we shall prove the main result of the
paper, Theorem
\ref{theorem:mainthm}.

\begin{acknowledgment} The author would like to thank Wei-Min Chen, Ron Fintushel, Wladek Lorek
and Tom Parker for helpful conversations. In particular, the author would like to thank Wei-Min
Chen for stimulating his interest in this subject. He would also like to thank Ron Fintushel,
Wladek Lorek and Tom Parker for their comments regarding the proofs of vanishing theorems for
Seiberg-Witten invariants, which comments led the author to revise the original presentation of
this paper.
\end{acknowledgment}

\setcounter{equation}{0}

\section{Restrictions on Connect Sum Decompositions of $M$}
\label{sec:restrictions}

In this section, we shall prove some technical lemmas about connect sum decompositions of
closed connected oriented $3$-manifolds whose product with a circle admits a symplectic
structure. Throughout this section, $M$ denotes a closed connected oriented
$3$-manifold and $X$ denotes the product $\mxs$ of $M$ with the circle $S^1$. We assume that
$S^1$ is equipped with the standard orientation, and $X$ is equipped with the corresponding
product orientation.

If $Y$ is a topological space and $j$ is a nonnegative integer, then
$H_j(Y)$ will denote the $j$th homology group of $Y$ with integer coefficients, $H^j(Y)$ will
denote the
$j$th cohomology group of
$Y$ with integer coefficients, and $b_j(Y)$ will denote the $j$th Betti number of $Y$, (i.e.
the rank of a maximal torsion free subgroup of $H_j(Y)$).

The intersection pairing $Q$ is an integer valued, symmetric bilinear form on $\htx$ which
descends to a unimodular form on the quotient of $\htx$ by its torsion subgroup.
$\bnx$ denotes the rank of a maximal subgroup of $\htx$ on which $Q$ is negative definite, and
$\bpx$ denotes the rank of a maximal subgroup of $\htx$ on which $Q$ is positive definite. By
definition, the signature $\s$ of $Q$ is equal to $\bpx - \bnx$. Note that
$\btx = \bnx + \bpx$.

\begin{lemma} Let $M$ be a closed connected oriented $3$-manifold and $X$ be the product
$\mxs$ of
$M$ with the circle $S^1$. Then $\bnx = \bpx = \bom$.
\label{lemma:bettirelns} \end{lemma}

\begin{proof} By Poincare duality, $\bom = \btm$. On the other hand, by the Kunneth formula,
$\htx$ is isomorphic to the direct sum of $\htmthzs$ and $\homthos$. Hence,
$\btx = \btm \times b_0(S^1) + \bom \times b_1(S^1) = \btm + \bom = 2\bom$.

Let $a$ and $b$ be elements of $\htx$ corresponding to the ``summand''
$\htmthzs$ of
$\htx$. We may represent $a$ and $b$ by immersed oriented surfaces $A$ and
$B$ in $X$ whose images lie in $M \times \{-1\}$ and $M \times \{1\}$ respectively. Since
these representatives,
$A$ and $B$, are disjoint, $Q(a,b) = 0$.

Let $c$ and $d$ be elements of $\htx$ corresponding to the ``summand''
$\homthos$ of
$\htx$. We may represent $c$ and $d$ by immersed surfaces $C \times S^1$ and $D \times S^1$,
where $C$ and $D$ are immersed oriented circles in $M$. We may assume that
$C$ and $D$ are disjoint, so that these representatives, $C \times S^1$ and $D \times S^1$, of
$c$ and $d$ are disjoint. Again, we conclude that $Q(c,d) = 0$.

Hence, the intersection pairing $Q$ is zero on each of the two ``summands''
$\htmthzs$ and
$\homthos$ of $\htx$. It follows that the signature $\s$ of $Q$ is equal to
$0$. On the other hand, $\s = \bpx - \bnx$. Hence, $\bnx = \bpx$. It follows that $2\bom =
\btx = \bnx + \bpx = 2\bpx$ and, hence, $\bpx = \bom$.

This completes the proof of Lemma \ref{lemma:bettirelns}.
\end{proof}

\begin{lemma} Let $M$ be a compact connected $3$-manifold. Suppose that
$\ab$ is a connect sum decomposition of $M$ such that $\boa \geq 1$ and $B$ has a nontrivial
connected covering space with finite index. Then there exists a connected covering space $\tm$
of
$M$ with finite index and an embedded $2$-sphere $\F$ in $\tm$ such that $\botm
\geq 2$ and $\F$ is nonseparating, (i.e. $\tmsS$ is connected).
\label{lemma:badcover} \end{lemma}

\begin{proof} Since $M$ is compact and connected, $A$ and $B$ are each compact and connected.
By assumption, there exists a connected covering space $\tb$ of $B$ of finite index $p \geq 2$.
Since $\boa \geq 1$ and $A$ is a compact $3$-manifold, $\hoa$ is a finitely generated abelian
group with a torsion free subgroup of rank at least $1$. By the classification of finitely
generated abelian groups, there exists an infinite cyclic quotient $\Z$ of
$\hoa$. Since $\Z_p$ is a quotient of $\Z$ and $\hoa$ is isomorphic to the abelianization of
the fundamental group
$\fund(A)$ of $A$, we obtain an epimorphism $\lambda_A : \fund(A)
\rightarrow \Z_p$. By covering space theory, $\lambda_A$ corresponds to a connected covering
space $\ta$ of $A$ with index $p$.

By assumption, $M$ is obtained from the disjoint union of $A$ and $B$ by removing a ball $\va$
from
$A$ and a ball $\vb$ from $B$ and gluing the resulting complements $A
\setminus \va$ and
$B \setminus \vb$ together along their boundaries $\bva$ and $\bvb$ by a diffeomorphism
$\phi : \bva \rightarrow \bvb$.

Since $\va$ is simply connected, the preimage $\tva$ of $\va$ in $\ta$ is a disjoint union of
$p$ balls, $\vaj, 1 \leq j \leq p$, in $\ta$. Likewise, the preimage $\tvb$ of
$\vb$ in $\tb$ is a disjoint union of $p$ balls, $\vbj, 1 \leq j \leq p$, in $\tb$. Note that
the diffeomorphism
$\phi : \bva \rightarrow \bvb$ lifts to a diffeomorphism $\phi_j : \bvaj
\rightarrow \bvbj$, for each $j$ with $1 \leq j \leq p$. Hence, we may form a compact connected
$3$-manifold $\tm$ by removing $\tva$ from $\ta$ and $\tvb$ from $\tb$, and gluing $\bvaj$ to
$\bvbj$ by the diffeomorphism $\phi_j : \bvaj \rightarrow \bvbj$, for each $j$ with $1
\leq j \leq p$. Note that
$\tm$ is a connected covering space of $M$ with index $p$.

Since $\ta$ is a connected $3$-manifold and $\tva$ is a disjoint union of
$p$ balls in
$\ta$, $\ta \setminus \tva$ is a connected $3$-manifold. Likewise,
$\tb \setminus \tvb$ is a connected $3$-manifold. Note that $\partial V_{A,1}$ determines a
smoothly embedded $2$-sphere $\F$ in $\tm$. The manifold $N$ obtained by cutting
$\tm$ along $\F$ may be constructed by removing $\tva$ from $\ta$ and
$\tvb$ from $\tb$, and gluing $\bvaj$ to $\bvbj$ by the diffeomorphism
$\phi_j : \bvaj \rightarrow \bvbj$, for each $j$ with $2 \leq j \leq p$. Since $p \geq 2$, we
conclude that $N$ is a connected $3$-manifold. Since $\tmsS$ is equal to the interior of $N$,
it follows that $\tmsS$ is connected.

Note that $\tva$ determines an embedding of a disjoint union of $p$
$2$-spheres in $\tm$. We may apply the Mayer-Vietoris sequence to the corresponding
decomposition of
$\tm$ into two submanifolds
$\ta \setminus \tva$ and $\tb \setminus \tvb$. In particular, we conclude that the first Betti
number $\botm$ of $\tm$ satisfies the equation $\botm = \bota + \botb + p - 1$. (Here, we
observe that, by the Mayer-Vietoris sequence, the removal of $p$ disjoint balls from $\ta$
($\tb$) does not affect $\hota$ ($\hotb$).) Since $\ta$ is a covering space of $A$ with finite
index, $\bota \geq \boa$. Indeed, the covering projection $\pi : \ta
\rightarrow A$ induces a homomorphism $\pi_* : \hota \rightarrow \hoa$ which maps $\hota$
onto a subgroup of finite index in $\hoa$. Thus, since $\boa \geq 1$, $\bota \geq 1$. Since
$\botb \geq 0$ and
$p \geq 2$, we conclude that
$\botm = \bota + \botb + p - 1 \geq 1 + 0 + 2 - 1 \geq 2$.

This completes the proof of Lemma \ref{lemma:badcover}.

\end{proof}

\begin{lemma} Let $M$ be a closed connected oriented $3$-manifold for which the product $X =
\mxs$ of $M$ with the circle $S^1$ admits a symplectic structure. Suppose that
$\ab$ is a connect sum decomposition of $M$ such that $\boa \geq 1$. Then every connected
covering space of $B$ with finite index is trivial.
\label{lemma:mainlemma} \end{lemma}

\begin{remark} The main facts used in the proof of Lemma
\ref{lemma:mainlemma} are familiar: (i) the vanishing of the Seiberg-Witten invariants for
$4$-manifolds admitting appropriate decompositions, and (ii) the nonvanishing of certain
Seiberg-Witten invariants for closed symplectic $4$-manifolds.
\label{lemma:famfacts} \end{remark}

\begin{remark} In order to apply the familiar facts mentioned in the previous remark, we must
pass to an appropriate covering space $\tm$ of $M$. The idea of passing to a covering space to
apply these familiar facts was introduced by Dieter Kotschick in \cite{k1}, and exploited
further in
\cite{kmt} and
\cite{k3}. Kotschick's ``covering trick'' exploits a particular covering space $\tx$ of a
$4$-manifold $X$. The covering trick which we shall use in the proof of Lemma
\ref{lemma:mainlemma} is similar to Kotschick's covering trick, but the covering space $\tm$ of
$M$ which we employ is not the three-dimensional analogue of the covering space $\tx$ of
$X$ exploited by Kotschick. The difference between our choice of covering spaces and
Kotschick's choice of covering spaces corresponds to the fact that we appeal to a different
vanishing result for Seiberg-Witten invariants than that which is invoked in Kotschick's
covering trick.
\label{remark:compcovtrick} \end{remark}

\begin{proof}[Proof of Lemma \ref{lemma:mainlemma}] Suppose, on the contrary, that $B$ has a
nontrivial connected covering space with finite index. Let $\tm$ and $\F$ be as in Lemma
\ref{lemma:badcover}. Orient $\tm$ by pulling back the orientation on
$M$ through the covering map
$\pi : \tm \rightarrow M$. Orient the embedded $2$-sphere $\F$ in
$\tm$. Since $\F$ is a smooth oriented nonseparating embedded $2$-sphere in the oriented
$3$-manifold $\tm$, we may choose a smooth embedded oriented circle $\g$ in
$\tm$ such that
$\F$ meets $\g$ exactly at one point $x$ in $\tm$, this point is a point of transverse
intersection of $\F$ with $\g$ in $\tm$, and the sign of intersection of
$\F$ with
$\g$ at $x$ is positive.

Let $\tx$ denote the product $\tmxs$ of $\tm$ with the circle $S^1$. The smoothly embedded
oriented
$2$-sphere $\F$ in $\tm$ determines a smoothly embedded oriented $2$-sphere
$\ts = \F \times \{1\}$ in the closed oriented $4$-manifold $\tx$. Since
$\ts$ lies in the hypersurface $\tm \times \{1\}$ in $\tx$, $\ts$ has square zero (i.e.
$\tilde Q(\brts,\brts) = 0$, where $\tilde Q$ is the intersection pairing on $\tx$, and
$\brts$ is the homology class in
$H_2(\tx)$ represented by $\ts$). The oriented circle $\g$ in $\tm$ determines an oriented
torus
$\tT = \g \times S^1$ in $\tx$. By our previous assumptions, $\ts$ meets
$\tT$ in exactly one point, $(x,1)$, this point is a point of transverse intersection
of
$\ts$ with $\tT$ in $\tx$, and the sign of intersection of $\ts$ with $\tT$ at this point is
positive. Hence, $\tilde Q(\brts,[\tT]) = 1$. It follows that $\ts$ is {\it essential}, (i.e
$\brts$ is a homology class of infinite order in $H_2(\tx)$).

By Lemma \ref{lemma:bettirelns} and Lemma \ref{lemma:badcover}, $b_2^+(\tx) = \bom \geq 2$.
Thus, $\tx$ is a closed oriented $4$-manifold with $b_2^+(\tx) > 1$, and
$\ts$ is an essential embedded sphere in $\tx$ of nonnegative self-intersection. Hence, by
Lemma 5.1 in
\cite{fs}, the Seiberg-Witten invariant $SW_{\tx}$ vanishes identically.

Since $\tm$ is a covering space of $M$, $\tx$ is a covering space of $X =
\mxs$. Since, by assumption, $X$ admits a symplectic structure, $\tx$ admits a symplectic
structure. Indeed, we may pull back any symplectic structure $\w$ on $X$ through the covering
map $\pi :
\tx \rightarrow X$ to obtain a symplectic structure $\pi^* \w$ on $\tx$. Let $\tw$ be a
symplectic structure on $\tx$. (By assumption, $\tw \wedge \tw$ gives the orientation of
$\tx$.) Then, by the Main Theorem of \cite{t}, the first Chern class of the associated almost
complex structure on $\tx$ has Seiberg-Witten invariant equal to $\pm 1$.

This is a contradiction. Hence, every connected covering space of $B$ with finite index is
trivial.

This completes the proof of Lemma \ref{lemma:mainlemma}.
\end{proof}

\section{The Main Result}
\label{sec:mainresult}

In this section, we shall prove the main result of this paper, Theorem
\ref{theorem:mainthm}. As in the previous section, $M$ denotes a closed connected oriented
$3$-manifold and $X$ denotes the product $\mxs$ of $M$ with the circle $S^1$. We assume that
$S^1$ is equipped with the standard orientation, and $X$ is equipped with the corresponding
product orientation.

A $3$-manifold $P$ is {\it non-trivial} if it is not homeomorphic to the
$3$-sphere $S^3$. We recall that a non-trivial $3$-manifold $P$ is {\it prime} if there is no
decomposition $P = M_1 \# M_2$ of $P$ as a connect sum with $M_1$ and $M_2$ non-trivial. In
\cite{m}, Milnor showed that each closed connected oriented $3$-manifold $P$ has a unique
decomposition as a connect sum of prime factors:

\begin{milnorthmone} Every nontrivial closed connected oriented
$3$-manifold $P$ is isomorphic to a connect sum $P_1 \# .... \# P_k$ of prime manifolds. The
summands $P_i$ are uniquely determined up to order and isomorphism.
\label{theorem:milnorthmone} \end{milnorthmone}

We recall that a $3$-manifold $P$ is {\it irreducible} if every (tame)
$2$-sphere in $P$ bounds a
$3$-ball in $P$. The relationship between primitivity and irreducibility for closed connected
oriented
$3$-manifolds may be summarized by the following results from \cite{m}:

\begin{milnorlemone} With the exception of manifolds isomorphic to $S^3$ or
$\sxs$, a closed connected oriented $3$-manifold is prime if and only if it is irreducible.
\label{lemma:milnorlemone} \end{milnorlemone}

\begin{milnorlemtwo} $\sxs$ is prime.
\label{lemma:milnorlemtwo} \end{milnorlemtwo}

Note that $S^3$ is irreducible, by a theorem of Alexander \cite{a}, but not prime, since $S^3$
is trivial. On the other hand, $\sxs$ is prime, by Lemma 1 (Milnor), but not irreducible,
since the $2$-sphere $S^2 \times \{1\}$ in $\sxs$ does not bound a
$3$-ball in $\sxs$.

We recall that a topological space $T$ is {\it aspherical} if all the higher homotopy groups
$\pii(T), i > 1,$ are zero. Note that neither $S^3$ nor $\sxs$ is aspherical, since
$\pi_3(S^3)$ and $\pi_2(\sxs)$ are both infinite cyclic. The relationship between
irreducibility and asphericity for closed connected oriented $3$-manifolds may be summarized
by the following result from \cite{m}:

\begin{milnorthmtwo} For every non-trivial closed connected oriented irreducible $3$-manifold
$P$, the second homotopy group $\pi_2(P)$ of $P$ is zero. If the Poincar\'{e} hypothesis is
true, then, conversely, every such manifold $P$ with $\pi_2(P) = 0$ is irreducible.
\label{theorem:milnorthmtwo} \end{milnorthmtwo}

Suppose now that $P$ is a closed connected oriented prime $3$-manifold with first Betti number
$\bop \geq 1$. Since $\bop \geq 1$, $P$ is not diffeomorphic to $S^3$. Hence, by Lemma 1
(Milnor), $P$ is either diffeomorphic to $\sxs$ or $P$ is irreducible. Since $\bop
\geq 1$, the fundamental group $\fund(P)$ of $P$ is infinite and, hence, the universal cover
$\tp$ of $P$ is not compact. In \cite{m}, Milnor observes that if, in addition to the
hypotheses on $P$ in Theorem 2 (Milnor), the universal covering space $\tp$ of $P$ is not
compact, then $\tp$ is contractible and all the higher homotopy groups
$\pii(P), i > 1,$ are zero. In other words, in this situation, $P$ is aspherical. Hence, we
have the following consequence of the above results from \cite{m}:

\begin{lemma} Suppose that $P$ is a closed connected oriented prime
$3$-manifold with first Betti number
$\bop \geq 1$. Then either $P$ is diffeomorphic to $\sxs$ or
$P$ is irreducible and aspherical.
\label{lemma:primefactors} \end{lemma}

\begin{proof}[Proof of Theorem \ref{theorem:mainthm}] By assumption, $X =
\mxs$ admits a symplectic structure, $\w$. We may assume that the closed $2$-form
$\w$ on $X$ represents an integral cohomology class $[\w]$ in the second cohomology
$H^2(X)$ of $X$. Since
$\w \cup \w$ is a positive closed $4$-form on the oriented $4$-manifold
$X$, the Poincare dual $e = PD([\w])$ is an element of $\htx$ with $Q(e,e) > 0$. It follows
that
$\bpx \geq 1$. Thus, by Lemma \ref{lemma:bettirelns}, $\bom \geq 1$ and, hence, $M$ is
nontrivial. Therefore, by Theorem 1 (Milnor), there exists a connect sum decomposition
$M = \mspo$ of $M$ into prime summands $M_i$, which are uniquely determined up to order and
isomorphism. Note that $M = S^3 \# \mspo$.

Since $M = \mspo$, $\bom = b_1(M_1) + ... + b_1(M_r)$. It follows, from the fact that
$\bom \geq 1$, that $b_1(M_i) \geq 1$ for some integer $i$ with  $1 \leq i
\leq r$. We may assume, without loss of generality, that $i = r$. Let $A = M_r$ and $B =
\snms$. Then
$M = \ab$, $A = M_r$ is prime, and $\boa = b_1(M_r) \geq 1$. By Lemma
\ref{lemma:primefactors}, either $A$ is diffeomorphic to $\sxs$ or $A$ is irreducible and
aspherical. By Lemma
\ref{lemma:mainlemma}, every connected covering space of $B$ with finite index is trivial.

This proves the existence of a connect sum decomposition of $M$ with the stipulated properties.

Suppose that $M = C \# D$ is a connect sum decomposition of $M$ where (i)
$b_1(C) \geq 1$, (ii) either $C$ is diffeomorphic to $\sxs$ or $C$ is irreducible and
aspherical, and (iii) every connected covering space of $D$ with finite index is trivial.

Since either $C$ is diffeomorphic to $\sxs$ or $C$ is irreducible and aspherical, it follows
from Lemma 1 (Milnor) and Lemma 2 (Milnor) that $C$ is prime. By Theorem 1 (Milnor), on the
other hand, there exists a connect sum decomposition $D = S^3 \# \ds$ of $D$ into prime
summands $D_i$. (Here, we allow for the possibility that $D$ is trivial, (i.e. $s = 0$).) It
follows that $M = C \# \ds$ is a decomposition of $M$ into prime summands, $C, D_1,..., D_s$.
By the uniqueness clause of Theorem 1 (Milnor), the factors $C, D_1,..., D_s$ must be
isomorphic to the factors $M_1,..., M_r$, respectively, up to a reordering of these factors.

By Lemma \ref{lemma:mainlemma}, every connected covering space of $B$ with finite index is
trivial. As pointed out in Remark \ref{remark:homologysphere}, this implies that $B$ is a
homology sphere. In particular, $b_1(B) = 0$. On the other hand, since $B = \snms$ and
$b_1(S^3) = 0$,
$b_1(B) = b_1(M_1) + ... + b_1(M_{r-1})$. It follows that $b_1(M_i) = 0$ for $1 \leq i < r$.
Thus,
$M_r$ is the unique prime factor of $M$ with positive first Betti number. Since $C$ is such a
prime factor of $M$, we conclude that $C = M_r = A$. Likewise, it follows that
$D = S^3 \# \ds = \snms = B$.

This proves the uniqueness of a connect sum decomposition of $M$ with the stipulated
properties.

This completes the proof of Theorem \ref{theorem:mainthm}.
\end{proof}


\bibliographystyle{amsalpha}

\end{document}